\documentclass{amsart}
\usepackage{amsfonts}
\usepackage{geometry}
\usepackage[onehalfspacing]{setspace}

\theoremstyle{plain}
\newtheorem{theorem}{Theorem}[section]
\newtheorem{lemma}[theorem]{Lemma}

\numberwithin{equation}{section}
\input{tcilatex}
\begin{document}
\title{\ \vspace*{0.25in}\\
Homological Properties of Color Lie Superalgebras}
\author{Kenneth L. Price}
\maketitle

\begin{abstract}
Let $\mathcal{L}=\mathcal{L}_{+}\oplus \mathcal{L}_{-}$ be a finite
dimensional color Lie superalgebra over a field of characteristic $0$ with
universal enveloping algebra $U(\mathcal{L})$. We show that $\limfunc{gldim}%
(U(\mathcal{L}_{+}))=\limfunc{lFPD}(U(\mathcal{L}))=\limfunc{rFPD}(U(%
\mathcal{L}))=\limfunc{injdim}_{U(\mathcal{L})}(U(\mathcal{L}))=\dim (%
\mathcal{L}_{+})$. We also prove that $U(\mathcal{L})$ is
Auslander-Gorenstein and Cohen-Macaulay and thus that it has a QF classical
quotient ring.
\end{abstract}

\setcounter{page}{287}

\section{Introduction}

Color Lie superalgebras are graded over an Abelian group $G$ and generalize
Lie superalgebras. Background information on Lie superalgebras and color Lie
superalgebras can be found in \cite{Kac} and \cite{Sc1} and \cite{BMPZ} and 
\cite{Sc2}, respectively. \ We will show, in a forthcoming paper \cite{Pr},
that the ideal structure of the enveloping algebra of a color Lie
superalgebra can be very different from enveloping algebras of Lie
superalgebras. \ By contrast, the results of this note illustrate that
certain homological properties are the same. \ Our main strategy is to pass
to the case where the grading group $G$ is finitely generated, so that the
color Lie superalgebra is determined by a Lie superalgebra and a 2-cocycle
defined on $G$.

In \S 2 we define color Lie superalgebras and state a theorem due to
Scheunert (see \cite{Sc2}). \ In \S 3 we calculate the finitistic and
injective dimensions of the enveloping algebra of a finite dimensional color
Lie superalgebra over a field of characteristic $0$. \ The enveloping
algebra of a finite dimensional color Lie superalgebra may have infinite
global dimension (this is known for ordinary Lie superalgebras, see \cite[%
Proposition 5]{Behr}). \ However, for a finite dimensional color Lie
superalgebra $\mathcal{L}$ which is \emph{positively graded}, i.e. $\mathcal{%
L}=\mathcal{L}_{+}$, we prove that, in analogy with the nongraded case, $%
\limfunc{gldim}U(\mathcal{L})=\dim \mathcal{L}$. \ In particular, theorem %
\ref{Homo Dimensions} generalizes \cite[Proposition 2.3]{KK}.

In theorem \ref{A-G C-M} we show that the enveloping algebra $U(\mathcal{L})$
of a color Lie superalgebra is Auslander-Gorenstein and Cohen-Macaulay. \ It
thus follows from \cite[Theorem 1.4]{KK2} that $U(\mathcal{L})$ has a (right
and left) quasi-Frobenius (QF) classical quotient ring.\vspace{0.2in}

\section{Color Lie Superalgebras}

Throughout $k$ denotes a field of characteristic $\neq 2$, $G$ denotes an
Abelian group, and all algebras are associative $k$-algebras with 1. We call
a map $\gamma :G\times G\rightarrow k^{\times }$ ($k^{\times }=k\backslash
\{0\}$) a \emph{skew-symmetric bicharacter} on $G$ if it satisfies\label%
{bicharacter 1} 
\begin{equation*}
\begin{array}{lll}
\text{(1)} & \text{ } & \gamma (f,gh)=\gamma (f,g)\gamma (f,h)\text{ and }%
\gamma (gh,f)=\gamma (g,f)\gamma (h,f)\text{ for any }f,g,h\in G\text{.} \\ 
\text{(2)} & \text{ } & \gamma (g,h)\gamma (h,g)=1\text{ for any }g,h\in G%
\text{.}%
\end{array}%
\end{equation*}

Note that (2) implies that $\gamma (g,g)=\pm 1$ for any $g\in G$. Set $%
G_{\pm }=\{g\in G:\,\gamma (g,g)=\pm 1\}$, then $G_{+}$ is a subgroup of $G$
and $[G:G_{+}]\leq 2$. \ In \cite[Lemma 2]{Sc2} it is shown that if $G\,$is
finitely generated then the skew-symmetric bicharacters on $G$ are
completely determined by the 2-cocycles on $G$ with coefficients in $%
k^{\times }$.

Generally, we shall be discussing objects which are graded over $G$ so let $G
$\textbf{-vec} denote the category of $G$-graded vector spaces and graded
linear maps. \ For any $V\in G$\textbf{-vec}, set $\partial x=g$ if $0\neq
x\in V_{g}$. \ Let $\gamma $ be a skew-symmetric bicharacter on $G$. \ We
will shorten our notation by writing $\gamma (x,y)$ instead of $\gamma
(\partial x,\partial y)$ for homogeneous $0\neq x\in V$, $0\neq y\in W$ and $%
V,W\in G$\textbf{-vec}. \ Note that $G$\textbf{-vec} is naturally contained
in $\mathbb{Z}_{2}$\textbf{-vec} via the decomposition $G=G_{+}\dot{\cup}%
G_{-}$ and the group homomorphism $\pi :G\rightarrow \mathbb{Z}_{2}$ given
by $\pi (G_{+})=\bar{0}$ and $\pi (G_{-})=\bar{1}$.\smallskip 

\noindent \textbf{Definition.~~} A $(G,\gamma )$\emph{-color Lie superalgebra%
} is a Lie algebra a pair $(\mathcal{L},\langle ,\rangle )$ such that $%
\mathcal{L}\in G$\textbf{-vec} and $\langle ,\rangle :\mathcal{L}\otimes _{k}%
\mathcal{L}\rightarrow \mathcal{L}$ is a graded bilinear map which satisfies
the following for any homogeneous $x,y,z\in \mathcal{L}$. 
\begin{equation*}
\begin{tabular}{ll}
$\gamma \text{\emph{-skew-symmetry}}$ & $\langle x,y\rangle =-\gamma
(x,y)\langle y,x\rangle $ \\ 
$\gamma \text{\emph{-Jacobi identity}}$ & $\gamma (z,x),\langle x,\langle
y,z\rangle \rangle +\gamma (y,z),\langle z,\langle x,y\rangle \rangle
+\gamma (x,y),\langle y,\langle z,x\rangle \rangle =0$%
\end{tabular}%
\end{equation*}

\noindent \textbf{Example. \ } (1) \ A \emph{Lie superalgebra} is a $(%
\mathbb{Z}_{2},\chi )$-color Lie superalgebra where $\chi (\overline{i},%
\overline{j})=$\newline
$\hspace*{0.47in}(-1)^{ij}$ for any $i,j\in \mathbb{Z}$.

\begin{enumerate}
\item[(2)] If $A$ is a $G$-graded algebra then there is a color Lie
superalgebra structure $A^{-}$\ defined on $A$ subject to the condition that 
$\langle a,b\rangle =ab-\gamma (a,b)ba$ for any nonzero homogeneous $a,b\in A
$.
\end{enumerate}

We shall now describe how to construct color Lie superalgebras from Lie
superalgebras as in \cite{Sc2}. \ In as much as Lie superalgebras are $%
\mathbb{Z}_{2}$-graded, we must consider $G$-gradings on Lie superalgebras
which respect the $\mathbb{Z}_{2}$-grading.\smallskip 

\noindent \textbf{Definition. \ \ }A Lie superalgebra $(L,[,])$, $L=L_{\bar{0%
}}\bigoplus L_{\bar{1}}$, is $G$\emph{-graded} if $L=\bigoplus_{g\in G}L_{g}$%
, $[L_{f},L_{g}]\subseteq L_{fg}$ for any $f,g\in G$, and $L_{h}\subseteq L_{%
\bar{0}}$ or $L_{h}\subseteq L_{\bar{1}}$ for each $h\in G$.

Starting with a Lie superalgebra $L$ there are many choices of groups $G$
and $G$-gradings on $L$. \ Let $G(L)$ be the smallest subgroup of $G$ which
grades $L$, that is, $G(L)$ is generated by $\left\{ g\in G:L_{g}\neq
0\right\} $. \ We cannot always pass to the case that $G=G(L)$ since, for
example, there are $G$-graded representations of $L$ which are not $G(L)$%
-graded.

Let $G\left( L\right) _{+}$ be the subgroup of $G(L)$ which is generated by $%
\{g\in G:U(L)_{g}\subseteq U(L)_{\bar{0}}\}$. \ Then $[G(L):G(L)_{+}]\leq 2$%
. \ Set $G\left( L\right) _{-}=G\left( L\right) \backslash G\left( L\right)
_{+}$; then the pair $(L,[,])$ is a $(G\left( L\right) ,\gamma _{0})$-color
Lie superalgebra where $\gamma _{0}$ is the skew-symmetric bicharacter on $%
G\left( L\right) $ defined below. 
\begin{equation*}
\gamma _{0}(g,h)=\left\{ 
\begin{array}{rl}
-1 & \text{if }g,h\in G\left( L\right) _{-} \\ 
1 & \text{otherwise}%
\end{array}%
\right. 
\end{equation*}%
To continue our construction of a color Lie superalgebra from a Lie
superalgebra, we need the notion of a 2-cocycle.\smallskip 

\noindent \textbf{Definition. \ \ }A \emph{2-cocycle} on $G$ is a map $%
\sigma :G\times G\rightarrow k^{\times }$ which satisfies 
\begin{equation*}
\sigma (f,gh)\sigma (g,h)=\sigma (f,g)\sigma (fg,h)
\end{equation*}%
for any $f,g,h\in G$.\label{cocycle}

If $\sigma $ is a 2-cocycle on $G$ then there is a skew-symmetric
bicharacter $\gamma $ defined on $G\left( L\right) $ by $\gamma (g,h)=\gamma
_{0}(g,h)\sigma (g,h)\sigma (h,g)^{-1}$\label{epsilon} for any $g,h\in
G\left( L\right) $. \ Moreover, there is a $(G\left( L\right) ,\gamma )$%
-color Lie superalgebra $(L^{\sigma },[,]^{\sigma })$ which has the same
vector space structure as $L$ but bracket $[,]^{\sigma }:L^{\sigma }\times
L^{\sigma }\rightarrow L^{\sigma }$ defined subject to the condition that $%
[x,y]^{\sigma }=\sigma (\partial x,\partial y)[x,y]$ for any homogeneous $%
x,y\in L$.

By setting $f=h=e$ in Definition \ref{cocycle} (where $e$ is the identity
element of $G$) we obtain $\sigma (e,e)=\sigma (g,e)=\sigma (e,g)$ for any $%
g\in G$. \ Note that there is no loss in generality by assuming that $\sigma
(e,e)=1$ since we can replace $\sigma $ with $\sigma ^{\prime }$, where $%
\sigma ^{\prime }=\sigma (e,e)^{-1}\sigma $ and have $L^{\sigma }\cong
L^{\sigma ^{\prime }}$ as $(G\left( L\right) ,\gamma )$-color Lie
superalgebras.

Theorem \ref{Scheunert's Theorem} summarizes the results from \cite{Sc2} we
are interested in. \ For the reader's convenience we will go over some basic
definitions.\smallskip 

\noindent \textbf{Definition. \ \ }Let $(\mathcal{L},\langle ,\rangle )$ be
a $(G,\gamma )$\emph{-color Lie superalgebra}.

\begin{enumerate}
\item A linear map $\phi :\mathcal{L}_{1}\rightarrow \mathcal{L}_{2}$
between $(G,\gamma )$-color Lie superalgebras $(\mathcal{L}_{1},\langle
,\rangle _{1})$ and $(\mathcal{L}_{2},\langle ,\rangle _{2})$ is called a 
\emph{homomorphism} if $\phi (\langle x,y\rangle _{1})=\langle \phi (x),\phi
(y)\rangle _{2}$ for any $x,y\in \mathcal{L}_{1}$.

\item A $G$\emph{-graded representation }of $\mathcal{L}$ is a pair $(V,\rho
)$ where $V\in G$\textbf{-vec} and $\rho :\mathcal{L}\rightarrow
End_{k}(V)^{-}$ is a homomorphism of color Lie superalgebras ($End_{k}(V)$
is a $G$-graded algebra since $V\in G$\textbf{-vec}$).$

\item For any $G$-graded algebra $A$, a graded map $\phi :\mathcal{L}%
\rightarrow A^{-}$ is called \emph{compatible} if it is homomorphism of $%
(G,\gamma )$-color Lie superalgebras.

\item The \emph{universal enveloping algebra} of $\mathcal{L}$ is a $G$%
-graded algebra $U(\mathcal{L})$ and a compatible map $\iota :\mathcal{L}%
\rightarrow U(\mathcal{L})$ which satisfies the property that for any
compatible map $\phi :\mathcal{L}\rightarrow A$ there is a unique (graded)
algebra homomorphism $\Phi :U(\mathcal{L})\rightarrow A$ such that $\Phi
\circ \iota =\phi $.
\end{enumerate}

\begin{theorem}[Scheunert]
Let $\mathcal{L}$ be a $(G,\gamma )$-color Lie superalgebra and $L$ and $%
L^{\sigma }$ as above.\label{Scheunert's Theorem}

\begin{enumerate}
\item If $G$ is finitely generated, then any color Lie superalgebra can be
obtained as an $L^{\sigma }$ for appropriately chosen $L$ and 2-cocycle $%
\sigma $.

\item The enveloping algebra $U(L^{\sigma })$ is obtained from $U(L)$ by
defining a new multiplication $\ast $ on $U$ subject to the condition that
for any homogeneous $x,y\in U(L)$ we have $x\ast y=\sigma (\partial
x,\partial y)xy.$

\item For a graded representation $\rho :L\rightarrow End_{k}(V)$ of $L$,
there is a graded representation $\rho ^{\sigma }:L^{\sigma }\rightarrow
End_{k}(V)$ of $L^{\sigma }$ which is obtained from $\rho $ subject to the
condition that for any homogeneous $x\in L^{\sigma }$ and $v\in V$, $\rho
^{\sigma }(x)(v)=\sigma (\partial x,\partial v)\rho (x)(v)$. \ This defines
a category equivalence between the categories of graded representations of $%
L $ and $L^{\sigma }$.
\end{enumerate}
\end{theorem}

We are particularly interested in parts (2) and (3) of theorem \ref%
{Scheunert's Theorem}. \ More generally, consider an arbitrary $G$-graded
algebra $\mathcal{R}$ and left $\mathcal{R}$-modules $V$ and $W$. \ An $%
\mathcal{R}$-module homomorphism $\phi :V\rightarrow W$ is called \emph{%
graded} if $\phi (V_{h})\subseteq W_{h}$ for each $h\in G$. \ Let $_{R}%
\mathcal{M}^{G}$ denote the category whose objects are $G$-graded left $%
\mathcal{R}$-modules and morphisms are all graded $\mathcal{R}$-module
homomorphisms from $V$ to $W$, denoted $HOM_{\mathcal{R}}^{0}(V,W)$, where $%
V,W\in _{R}\mathcal{M}^{G}$. \ In particular, $_{k}\mathcal{M}^{G}=G$\textbf{%
-vec}.

\begin{lemma}
Let $\sigma $ be a 2-cocycle on $G$ and $R$ a $G$-graded $k$-algebra.\label%
{cat equiv}

\begin{enumerate}
\item There is a $G$-graded $k$-algebra $R^{\sigma }$ with the same vector
space structure as $R$ and whose multiplication \thinspace $*$ is obtained
from the multiplication of $R$ subject to the condition $r*s=\sigma
(\partial r,\partial s)rs$ for any homogeneous $r,s\in R$.

\item For any $M\in _{R}\mathcal{M}^{G}$ there corresponds $M^{\sigma }\in
_{R^{\sigma }}\mathcal{M}^{G}$ such that $M^{\sigma }$ has the same vector
space structure as $M$ but the $R^{\sigma }$-module structure $.$ is
obtained from the $R$-module structure of $M$ subject to the condition $%
r.m=\sigma (\partial r,\partial m)rm$ for any homogeneous $r\in R$ and $m\in
M$.

\item The functor $^{\sigma }:_{R}\mathcal{M}^{G}\rightarrow _{R^{\sigma }}%
\mathcal{M}^{G}$ defined as in (2) is a category equivalence. In particular, 
$HOM_{R}^{0}(V,W)=HOM_{R^{\sigma }}^{0}(V^{\sigma },W^{\sigma })$. \
Moreover, if $V\in _{R}\mathcal{M}^{G}$ is graded free, then so is $%
V^{\sigma }\in _{R^{\sigma }}\mathcal{M}^{G}$.

\item For any $V,W\in _{R}\mathcal{M}^{G}$ and $\phi \in HOM_{R}^{0}(V,W)$
we have $\ker (\phi ^{\sigma })=(\ker \phi )^{\sigma }=\ker \phi $ and $%
\limfunc{im}(\phi ^{\sigma })=(\limfunc{im}\phi )^{\sigma }=\limfunc{im}\phi 
$.
\end{enumerate}
\end{lemma}

\proof%
It is easy to prove (1)-(3) directly. For (4), first note that $\ker \phi $, 
$\limfunc{im}\phi $, $\ker (\phi ^{\sigma })$ and $\limfunc{im}(\phi
^{\sigma })$ are $G$-graded submodules. \ For any homogeneous $x\in V$ we
have $x\in \ker (\phi ^{\sigma })\Leftrightarrow x\in \ker (\phi
)\Leftrightarrow x\in (\ker \phi )^{\sigma }$. \ This proves that $\ker
(\phi ^{\sigma })=(\ker \phi )^{\sigma }=\ker \phi $. \ The proof that $%
\limfunc{im}(\phi ^{\sigma })=(\limfunc{im}\phi )^{\sigma }=\limfunc{im}\phi 
$ is similar.%
\endproof%

\section{Homological Properties of $U(\mathcal{L})$}

Throughout this section, $\gamma $ is a skew-symmetric bicharacter on $G$
and $G_{+}$=$\left\{ g\in G:\gamma \left( g,g\right) =\pm 1\right\} $. For
any $V\in G$\textbf{-vec}, set $V_{\pm }=\oplus _{g\in G_{\pm }}V_{g}$. If $%
V=V_{+}$ (respectively $V=V_{-}$) then $V$ is called \emph{positively }%
(respectively\emph{\ negatively})\emph{\ graded}. We shall be primarily
interested in the case where $\mathcal{L}$ is an Abelian color Lie
superalgebra. \ If $\mathcal{L}$ is a positively graded, finite dimensional
and Abelian $(G,\gamma )$-color Lie superalgebra, then $U(\mathcal{L})\cong
k^{\gamma }[\theta _{1},\theta _{2},\ldots ,\theta _{n}]$, the \emph{color
polynomial ring} in $n=\dim \mathcal{L}$ variables (see \cite{BMPZ}). \ If $%
\mathcal{L}$ is a negatively graded, finite dimensional and Abelian $%
(G,\gamma )$-color Lie superalgebra, then $U(\mathcal{L})\cong \Lambda
^{\gamma }(\mathcal{L})$, the \emph{color exterior algebra (or color
Grassmann algebra) of }$\mathcal{L}$ (see \cite{BMPZ}).

\noindent \textbf{Example. }

\begin{enumerate}
\item Choose $q\in k^{\times }$. Let $\gamma $ be the skew-symmetric
bicharacter on $G=\mathbb{Z}^{2}$ defined by $\gamma
(g,h)=q^{g_{1}h_{2}-g_{2}h_{1}}$ for any $g=(g_{1},g_{2})$, $%
h=(h_{1},h_{2})\in G$.  Then $G=G_{+}$. \ Consider the Abelian $(G,\gamma )$%
-color Lie superalgebra with homogeneous basis $\{x,y\}$ such that $\partial
x=(1,0)$ and $\partial y=(0,1)$. \ Then $\mathcal{L}$ is positively graded
and $U(\mathcal{L})\cong k\left[ x,y:xy=qyx\right] $, the so-called \emph{%
quantum plane}.

\item Choose $0\neq q\in k^{\times }$. \ Let $\gamma $ be the skew-symmetric
bicharacter on $G=\mathbb{Z}^{2}$ defined by $\gamma
(g,h)=(-1)^{g_{1}h_{1}+g_{2}h_{2}}q^{g_{1}h_{2}-g_{2}h_{1}}$ for any $%
g=(g_{1},g_{2})$, $h=(h_{1},h_{2})\in G$. \ Then $G_{+}=\{(i,j)\in G:i+j\in 2%
\mathbb{Z}\}$ and $G_{-}=\{(i,j)\in G:i+j-1\in 2\mathbb{Z}\}$. \ Consider
the Abelian $(G,\gamma )$-color Lie superalgebra with homogeneous basis $%
\{x,y\}$ such that $\partial x=(1,0)$ and $\partial y=(0,1)$. \ Then $%
\mathcal{L}$ is negatively graded and $U(\mathcal{L})\cong k[x,y:xy=qyx]/I$,
where $I$ is the ideal generated by $\{x^{2},y^{2}\}$.
\end{enumerate}

If $V,W\in G$\textbf{-vec} then the tensor product $V\bigotimes W\in G$%
\textbf{-vec} has grading defined by $(V\otimes W)_{g}=\oplus _{h\in
G}(V_{h}\otimes W_{gh^{-1}})$. For $G$-graded algebras $A$ and $B$ there is
a $\gamma $\emph{-graded tensor product} $A\widehat{\otimes }_{k}B$ with
multiplication defined subject to the condition that $(a\widehat{\otimes }%
b)(a^{\prime }\widehat{\otimes }b^{\prime })=\gamma (b,a^{\prime
})aa^{\prime }\widehat{\otimes }bb^{\prime }$ for any nonzero homogeneous $%
a,a^{\prime }\in A$ and $b,b^{\prime }\in B$. \ The $\gamma $-graded tensor
product is essential to the study of enveloping algebras of color Lie
superalgebras. \ For example, the enveloping algebra of a color Lie
superalgebra is not a Hopf algebra in the usual sense; it is only a Hopf
algebra with respect to the $\gamma $-graded tensor product (see \cite[\S %
3.2.9]{BMPZ} for details).

For any Noetherian ring $R$, the left (respectively right) finitistic
dimension of $R$, denoted $\limfunc{lFPD}\left( R\right) $ (respectively $%
\limfunc{rFPD}\left( R\right) $), is the supremum of the projective
dimensions of left (respectively right) $R$-modules of finite projective
dimension. The following theorem generalizes \cite[Proposition 2.3]{KK}. \
As in \cite{KK} and \cite{KK2}, we assume that $\limfunc{char}k=0$ in all
that follows.

\begin{theorem}
\label{Homo Dimensions}Let $\mathcal{L}=\mathcal{L}_{+}\oplus \mathcal{L}_{-}
$ be a finite dimensional color Lie superalgebra. Then $\limfunc{gldim}(U(%
\mathcal{L}_{+}))=\limfunc{lFPD}(U(\mathcal{L}))=\limfunc{rFPD}(U(\mathcal{L}%
))=\limfunc{injdim}_{U(\mathcal{L})}(U(\mathcal{L}))=\dim (\mathcal{L}_{+})$.%
\label{dimension theorem}
\end{theorem}

\proof%
Let $U=U(\mathcal{L})$ denote the universal enveloping algebra of $\mathcal{L%
}$ and set $n=\dim (\mathcal{L}_{+})$. \ We first show that $\limfunc{lFPD}%
(U)=\limfunc{rFPD}(U))=\limfunc{injdim}_{U}(U)=n\geq \limfunc{gldim}(U(%
\mathcal{L}_{+}))$. \ Then we show that $\limfunc{gldim}(U(\mathcal{L}%
_{+}))\geq n$, that is, there exists a $U(\mathcal{L}_{+})$-module with
projective dimension at least $n$. \ As in the case of ordinary Lie
algebras, the existence of such a module will be demonstrated using the 
\emph{Chevalley-Eilenberg complex} (see Chapter 7 of \cite{We} for
background).

\begin{enumerate}
\item[1)] With respect to the standard filtration, $gr\left( U\right) $ is
Noetherian.
\end{enumerate}

The filtration is defined by $U^{-1}=\{0\}$, $U^{0}=k$ and, for $m>0$, $%
U^{m} $ is spanned by all monomials of length $\leq m$. \ By \cite[%
Proposition III.2.8]{BMPZ} $U$ is Noetherian and the associated graded
algebra $gr\left( U\right) $ is isomorphic to $k^{\gamma }[\theta
_{1},\theta _{2},\ldots ,\theta _{n}]\widehat{\otimes }\Lambda ^{\gamma }(%
\mathcal{L}_{-})$, where $\Lambda ^{\gamma }(\mathcal{L}_{-})$ is the color
exterior algebra determined by the vector space $\mathcal{L}_{-}$ with an
Abelian color Lie super algebra structure. \ Thus $gr\left( U\right) $ is an
iterated Ore extension of $\Lambda =\Lambda ^{\gamma }(\mathcal{L}_{-})$
without derivations, that is, $gr\left( U\right) \cong \Lambda \lbrack
\theta _{1};\alpha _{\partial \theta _{1}}][\theta _{2};\alpha _{\partial
\theta _{2}}]\cdots \lbrack \theta _{n};\alpha _{\partial \theta _{n}}]$. \
Therefore, $gr\left( U\right) $ is Noetherian by \cite[Theorem I.2.10]{McR}.

\begin{enumerate}
\item[2)] $\limfunc{injdim}_{U}U\leq \limfunc{injdim}_{gr\left( U\right)
}gr\left( U\right) =n$ and $\limfunc{gldim}U(\mathcal{L}_{+})\leq \limfunc{%
gldim}gr(U(\mathcal{L}_{+}))=n$
\end{enumerate}

Since $U^{m}=0$ when $m<0$, the filtration $\{U^{m}\}$ satisfies the \emph{%
closure condition} of \cite[\S 1.2.11]{Bj}. \ Therefore, we may apply \cite[%
Theorem 1.3.12]{Bj} to obtain $\limfunc{injdim}_{U}U\leq \limfunc{injdim}%
_{gr\left( U\right) }gr\left( U\right) $. \ Moreover, $\limfunc{gldim}U(%
\mathcal{L}_{+})\leq \limfunc{gldim}gr(U(\mathcal{L}_{+}))$ by \cite[7.6.18]%
{McR}. For any ring $R$ and automorphism $\sigma $ of $R$, one can show that 
$\limfunc{gldim}(R[t;\sigma ])=\limfunc{gldim}(R)+1$ and $\limfunc{injdim}%
_{R[t;\sigma ]}(R[t;\sigma ])=\limfunc{injdim}_{R}(R)+1$ by similar methods
to the proofs in the case that $\sigma $ is the identity on $R$. \ Since $%
gr(U(\mathcal{L}_{+})\cong k^{\gamma }[\theta _{1},\theta _{2},\ldots
,\theta _{n}]$ we have $\limfunc{gldim}gr(U(\mathcal{L}_{+}))=n=\dim 
\mathcal{L}_{+}$. \ By \cite[Corollary 6.3]{FMS} $\Lambda $ is a Frobenius
algebra; hence $\limfunc{injdim}_{gr(U)}gr(U)=n=\dim \mathcal{L}_{+}$ as $%
gr(U)\cong \Lambda ^{\gamma }[\theta _{1},\theta _{2},\ldots ,\theta _{n}]$.

\begin{enumerate}
\item[3)] $\limfunc{lFPD}(U)=\limfunc{rFPD}(U))=\limfunc{injdim}_{U}(U)\geq 
\limfunc{gldim}(U(\mathcal{L}_{+}))$
\end{enumerate}

Since $\limfunc{gldim}U(\mathcal{L}_{+})$ is finite and $U$ is free over $U(%
\mathcal{L}_{+})$ on both sides, we may apply \cite[Theorem 1.4]{Co} to
obtain $\limfunc{gldim}U(\mathcal{L}_{+})\leq \limfunc{lFPD}U$. \ By \cite[%
Proposition 2.1]{KK} we have $\limfunc{injdim}_{U}U=\limfunc{lFPD}(U)=%
\limfunc{rFPD}(U)$.

\emph{The rest of the proof is dedicated to showing that if }$\mathcal{L}$%
\emph{\ is a positively graded finite-dimensional }$(G,\gamma )$\emph{-color
Lie superalgebra then }$n=\dim \mathcal{L}\leq \limfunc{gldim}U(\mathcal{L})$%
\emph{. }

\begin{enumerate}
\item[4)] Pass to the case that $\mathcal{L}=L^{\sigma }$ for some Lie
algebra $L$.
\end{enumerate}

The modules we will consider are graded over the subgroup $G(\mathcal{L})$
of $G$ generated by $\{g\in G:\mathcal{L}_{g}\neq 0\}$. Therefore, we can
pass to the case that $G$ is finitely generated by assuming that $G=G(%
\mathcal{L})$. \ By theorem \ref{Scheunert's Theorem}(1), there is a $G$%
-graded Lie algebra $L$ with bracket $[,]$ and a 2-cocycle $\sigma $ on $G$
such that $\mathcal{L}=L^{\sigma }$.

\begin{enumerate}
\item[5)] The Chevalley-Eilenberg complex $V_{\ast }(L)^{\underrightarrow{%
\varepsilon }}k$ is a complex in $_{U\left( L\right) }\mathcal{M}^{G}$.
\end{enumerate}

We must show that $V_{\ast }(L)^{\underrightarrow{\varepsilon }}k$ is a 
\emph{graded complex}, i.e., the modules are graded and all differentiations
are graded maps. Let $\Lambda L$ denote the exterior algebra of $L$ and $%
\Lambda ^{i}L$ denote the $i^{\text{th}}$ exterior algebra of $L$. Then $%
\Lambda L=\oplus _{i=1}^{n}\Lambda ^{i}L$ can be $G$-graded so that the
homogeneous component of degree $g\in G$ is $(\Lambda
L)_{g}=\sum_{i=1}^{n}(\Lambda ^{i}L)_{g}$ and $(\Lambda
^{i}L)_{g}=\sum_{g_{1}g_{2}\cdots g_{i}=g}L_{g_{1}}\wedge L_{g_{2}}\wedge
\cdots \wedge L_{g_{i}}$ where the sum is indexed by all appropriate $%
g_{1},g_{2},\ldots ,g_{i}\in G$.

We have $V_{i}(L)=U(L)\otimes _{k}\Lambda ^{i}L$ which is a graded free left 
$U(L)$-module such that $V_{i}(L)_{g}=\sum_{h\in G}(U\left( L\right)
_{h})\otimes _{k}(\Lambda ^{i}L)_{hg^{-1}}$ so we only need to check that
the \emph{augmentation map} $\varepsilon :V_{0}(L)\rightarrow k$ is graded
and, for $1\leq i\leq n$, $d_{i}:V_{i}(L)\rightarrow V_{i-1}(L)$ is graded.
The augmentation map $\varepsilon :V_{0}(L)=U\left( L\right) \rightarrow k$
is graded since it is induced from the Lie algebra map $L\rightarrow k$
which sends everything to zero, a graded map. It is also easy to see that $%
d:V_{1}(L)\rightarrow V_{0}(L)$ is graded since it is just the product map $%
u\otimes x\longmapsto ux$. \ For $1<i\leq n$, the map $d_{i}:V_{i}(L)%
\rightarrow V_{i-1}(L)$ is defined by the formula $d_{i}(u\otimes
x_{1}\wedge x_{2}\wedge \cdots \wedge x_{i})=\theta _{1}+\theta _{2}$ where $%
\theta _{1}$ and $\theta _{2}$ are defined below. 
\begin{eqnarray*}
\theta _{1} &=&\sum_{j=1}^{i}(-1)^{j+1}ux_{j}\otimes x_{1}\wedge x_{2}\wedge
\cdots \wedge \hat{x}_{j}\wedge \cdots \wedge x_{i} \\
\theta _{2} &=&\sum_{l<m}(-1)^{l+m}u\otimes \lbrack x_{l},x_{m}]\wedge
x_{1}\wedge x_{2}\wedge \cdots \wedge \hat{x}_{l}\wedge \cdots \wedge \hat{x}%
_{m}\wedge \cdots \wedge x_{i}
\end{eqnarray*}%
Therefore $d_{i}$ is graded since if $u,x_{1},x_{2},\ldots ,x_{i}$ are
homogeneous elements then $u\otimes x_{1}\wedge x_{2}\wedge \cdots \wedge
x_{i},\theta _{1}$ and $\theta _{2}$ are all homogeneous of the same degree.

\begin{enumerate}
\item[6)] The Chevalley-Eilenberg complex $V_{\ast }(\mathcal{L})^{%
\underrightarrow{\varepsilon }}k$ is a complex in $_{U}\mathcal{M}^{G}$.
\end{enumerate}

By parts (2), (3) and (4) of Lemma \ref{cat equiv}, there is a projective
resolution $V_{\ast }(\mathcal{L})$ of $k=k^{\sigma }$ over $U(\mathcal{L})$
such that $V^{i}(\mathcal{L})=V^{i}(L)^{\sigma }$. \ For $1\leq i\leq n$ the
module $V^{i}(L)$ is a graded free $U(L)$-module hence the module $%
V^{i}(L)^{\sigma }$ is a graded free $U(L^{\sigma })$-module by Lemma \ref%
{cat equiv}(3).

\begin{enumerate}
\item[7)] $\limfunc{gldim}U(\mathcal{L})\geq n$
\end{enumerate}

As in \cite[Exercise 7.7.2]{We}, there is a graded representation $\rho
:L\rightarrow End_{k}(\Lambda ^{n}L)$ defined as below for $0\neq y\in L$
and $x_{1},x_{2},\ldots ,x_{n}$ a (homogeneous) basis of $L$.

\begin{equation*}
\rho (y)(x_{1}\wedge x_{2}\wedge \cdots \wedge
x_{n})=\sum_{i=1}^{n}x_{1}\wedge x_{2}\wedge \cdots \wedge \lbrack
y,x_{i}]\wedge \cdots \wedge x_{n}
\end{equation*}

Let $M=\Lambda ^{n}L$ be the (graded) $U(L)$-module defined by the above
action. \ Then according to \cite{We}, $\func{Ext}_{U(L)}^{n}(k,M)\cong \ker
(d_{n}^{\ast })/\limfunc{im}(d_{n-1}^{\ast })\cong k$. \ By parts (2) and
(4) of Lemma \ref{cat equiv} and 6) above, we have $\ker d_{n}^{\ast }=(\ker
d_{n}^{\ast })^{\sigma }=\ker (d_{n}^{\ast })^{\sigma }$ and $\limfunc{im}%
d_{n-1}^{\ast }=(\limfunc{im}d_{n-1}^{\ast })^{\sigma }=\limfunc{im}%
(d_{n-1}^{\ast })^{\sigma }$. \ This implies that 
\begin{equation*}
\func{Ext}_{U}^{n}(k,M^{\sigma })\cong \left( \func{Ext}_{U\left( L\right)
}^{n}(k,M)\right) ^{\sigma }\cong k
\end{equation*}%
therefore $pd_{U}(M^{\sigma })\geq n$ hence $\limfunc{gldim}U\geq n$ as
desired.

The result follows from 3) and 7) above. 
\endproof%

\begin{theorem}
\label{A-G C-M}Let $\mathcal{L}$ be a finite dimensional color Lie
superalgebra. \ Then $U(\mathcal{L})$ is Auslander-Gorenstein and
Cohen-Macaulay and thus has a quasi-Frobenius classical quotient ring.
\end{theorem}

\proof%
We have $\Lambda $ is Auslander-Gorenstein and Cohen-Macaulay since $\Lambda 
$ is a Frobenius algebra (see \cite[Corollary 6.3]{FMS}). \ As gr$(U(L))$ is
an iterated Ore extension of $\Lambda $ where each iteration is of the form $%
R[x;\sigma ]$, it follows from \cite[Lemma (ii), p. 184]{LS} and \cite[%
Theorem 4.2]{E} that gr$(U(L))$ is Auslander-Gorenstein and Cohen-Macaulay.
\ We thus conclude that $U(\mathcal{L})$ is Auslander-Gorenstein and
Cohen-Macaulay by \cite[Theorem 1.4.1]{Bj} and \cite[Lemma 4.4]{SZ} (note
that Theorem \ref{Homo Dimensions} is needed so that $\limfunc{injdim}_{U(%
\mathcal{L})}U(\mathcal{L})<\infty $). \ The last statement follows from 
\cite[Theorem 1.4]{KK2}.%
\endproof%

\noindent \textbf{Acknowledgements}

The author is indebted to Jim Kuzmanovich for many helpful conversations
concerning this material.

\end{document}